\documentclass[12pt]{amsart}
\usepackage{amsfonts}
\usepackage[pdftex]{graphicx}

\newcommand{\Q}{{\mathbb Q}}

\newcommand{\C}{{\mathbb C}}

\renewcommand{\epsilon}{\varepsilon}
\renewcommand{\phi}{\varphi}

\newtheorem{Lemma}{Lemma}[section]
\newtheorem{Theorem}{Theorem}[section]
\newtheorem{Proposition}{Proposition}[section]

\newtheorem{Corollary}{Corollary}[section]
\newtheorem{Definition}{Definition}[section]
\newtheorem{Example}{Example}[section]
\newtheorem{Remark}{Remark}[section]

\begin{document}

\address{Department of Mathematics, University of Pittsburgh, 301 Thackeray Hall, Pittsburgh, PA 15260, USA.}
\author{Alexander Borisov}
\thanks{The research of the author was supported in part by the NSA grants H98230-08-1-0129, H98230-06-1-0034 and H98230-11-1-0148.}
\title[unramified planar self-maps]{On resolution of compactifications of unramified planar self-maps}

\begin{abstract} The goal of this paper is to approach the two-dimensional Jacobian Conjecture using ideas of birational algebraic geometry. We study the resolution of  rational self-map of the projective plane that comes from a hypothetical counterexample to the two-dimensional Jacobian Conjecture and establish several strong restrictions on its structure. In particular, we get a very detailed description of its Stein factorization. We also establish some combinatorial results on determinants of the Gram matrix of weighted trees and forests and apply them to study exceptional divisorial valuations of the field of rational function in two vairables with centers outside of the affine plane. We use some techniques of modern birational geometry, in particular adjunction inequalities, to further restrict the structure of possible counterexamples. Ultimately, we hope that this paper will pave the way for settling the two-dimensional Jacobian Conjecture using the techniques of modern birational geometry. 
\end{abstract}

\email{borisov@pitt.edu}

\maketitle

\section{Introduction}
The goal of this paper is to develop a new geometric approach to the two-dimensional Jacobian Conjecture. Some of the results may also be used to understand better the structure of the polynomial automorphisms of the plane. All varieties in this paper are over complex numbers.

This paper is written primarily for the algebraic geometers, by an algebraic geometer. So we take for granted the properties of canonical divisors on normal surfaces, while working out carefully some elementary combinatorial results.

Suppose $f(x,y)$ and  $g(x,y)$ are two polynomials with complex coefficients. The classical Jacobian Conjecture (due to Keller) asserts the following.

{\bf Conjecture.} (Jacobian Conjecture in dimension two) If the Jacobian of the pair $(f,g)$ is a non-zero constant, then the map $(x,y)\mapsto (f(x,y),g(x,y))$ is invertible. Note that the opposite is clearly true, because the Jacobian of any polynomial map is a polynomial, and, when the map is invertible, it must have no zeroes, so it is a constant.

The Jacobian Conjecture and its generalizations received considerable attention in the past, see for example \cite{AbhyankarTataLectures}, \cite{Essen}. It is notorious for its subtlety, having produced a substantial number of incorrect proofs by respectable mathematicians. 

From the point of view of a birational geometer, the most natural approach to the two-dimensional Jacobian Conjecture is the following. Suppose a counterexample exists. It gives a rational map from $P^2$ to $P^2.$ After a sequence of blow-ups of points, we can get a surface $X$ with two maps: $\pi : X \to P^2$ (projection onto the origin $P^2$) and $\phi : X \to P^2$ (the lift of the original rational map).

Note that $X$ contains a Zariski open subset isomorphic to $A^2$ and its complement, $\pi ^* ((\infty))$, is a tree of smooth rational curves. We will call these curves exceptional, or curves at infinity. The structure of this tree is easy to understand inductively, as it is built from a single curve $(\infty)$ on $P^2$ by a sequence of two operations: blowing up a point on one of the curves or blowing up a point of intersection of two curves. However, a non-inductive description is probably impossible, which is the first difficulty in this approach. Another difficulty comes from the fact that the exceptional curves on $X$ may behave very differently with respect to the map $\phi$. More precisely, there are four types of curves $E$.

type 1) $\phi (E) = (\infty )$
 
type 2) $\phi (E)$ is a point on $(\infty )$

type 3) $\phi (E)$ is a curve, different from $(\infty)$ 

type 4) $\phi (E)$ is a point not on $(\infty)$

From a first glance, the situation appears almost hopeless. One of the goals of this paper is to show that it is really not that bad. In particular, for a given graph of curves, one can essentially always tell which curves are of which type, and there is a fairly restrictive family of graphs that can potentially appear in a counterexample to the Jacobian Conjecture. Our main tools are the basic tools of algebraic geometry of surfaces: the intersection pairing and the adjunction formula. We are guided by some of the ideas of the log Minimal Model Program, but most of the proofs are relatively elementary and self-contained.

This paper is organized as follows. Section 2 is devoted to preliminary results on the graph of the exceptional curves. In section 3 we use slightly more subtle arguments to get further restrictions on that graph and the structure of the Stein factorization of the morphism $\phi$.  Section 4 is devoted to adjunction formulas and inequalities that generalize the arguments of section 3 to more general compactifications of the target plane. In section  5 we develop and use the theory of  the determinants of weighted trees and forests and get some further restrictions on the structure of the exceptional graph. In section 6 we try to attack the problem by  exploring  different compactifications of the target plane. Finally, section 7 contains some informal discussion of possible further developments of our approach.

{\bf Acknowledgments.} This paper is dedicated to the memory of V.A. Iskovskikh, who introduced the author to the beauty of birational geometry. The author is also indebted to David Wright and Ed Formanek for stimulating discussions related to the Jacobian Conjecture. 

\section{Preliminary Observations and Definitions}
We change the notation slightly from the Introduction.

Suppose $X= P^2$, $Y=P^2$ and $\phi _Y^X: X ---> Y $ is a rational map. Suppose further that on an open subset $A^2\subset P^2 =X$ the map $\phi_Y^X$ is defined, unramified, and $\phi_Y^X(A^2) \subseteq A^2 \subset P^2 =Y.$ By a sequence of blow-ups at smooth points, we get a surface $Z$ with a birational map $\pi : Z \to X $ and a generically finite map $ \phi^{Z}_Y : Z \to Y$ such that $\phi ^{Z}_Y = \phi ^X_Y \circ \pi .$ We will denote $\phi = \phi ^{Z}_Y .$
 
The blow-ups that lead to $Z$ can be done outside of $A^2\subset X.$ So $Z = A^2 \cup (\cup E_i)$, where $E_i$ are rational curves. The following proposition collects some straightforward observations. 

\begin{Proposition}

1) The curves $E_i$ form  a tree.

2) One of $E_i$ is $\pi ^{-1} (\infty),$ all others are mapped to points by $\pi $.

3) The classes of $E_i$ form a basis of the Picard group of $Z$.

\end{Proposition}

The structure of $Z$ is to a large extent determined by the graph of intersections of $E_i$. The vertices of this graph correspond to $E_i$-s and are usually labeled by $E_i^2.$ The edges correspond to the points of intersections of two different $E_i$-s. The graph is a tree.

This graph is not so easy to deal with because blowing up a point changes the self-intersections of the curves passing through it. Inspired by the Minimal Model Program, we consider a different labeling of this graph. We consider the augmented canonical class of $Z,$ $\bar{K}_{Z} = K_{Z} +\sum _i E_i$. It can be uniquely written as a linear combination of $E_i,$ $\bar{K}_{Z}= \sum _i a_iE_i$. We label the vertices of the intersection graph by these numbers $a_i.$

With this labeling we now describe what happens when a point is blown up, in any of the  intermediate steps in getting from $X$ to $Z.$ 

\begin{Proposition}
When a point is blown up, going from $Z'$ to $Z'',$ one of the following two operations is performed to the graph of the exceptional curves:

1) A new vertex is added to the graph, connected to one of the vertices. It is labeled $a_i+1,$ where $a_i$ is the label of the vertex it is connected to.

2) A new vertex is introduced on the edge connecting two vertices, ``breaking" the edge into two edges. The new vertex gets labeled with $a_i+a_j,$ where $a_i$ and $a_j$ are the labels of the two vertices it is connected to.
\end{Proposition}

{\bf Proof.} The first case corresponds to blowing up a point on one of the curves. The second case corresponds to blowing up an intersection of two curves. The augmented canonical class calculations are straightforward and are left to the reader. $\square$

Notice that once a vertex is created, its label never changes, which is in sharp contrast with the traditional labeling.

The following observation is true for any $Z,$ unrelated to the map $\phi .$  It is easily proven by induction on the number of exceptional curves, using the above proposition.

\begin{Proposition} For any two adjacent vertices $E_i,$ $E_j$ of the graph of $Z,$ $gcd(a_i,a_j)=1.$ In particular, no two adjacent vertices have even labels.
\end{Proposition}

The following example serves two purposes. It shows how the graph of $Z$ is constructed from the graph of $X=P^2,$ and we will use it to compare our labeling with the traditional self-intersection labeling.

{\bf Example.} 

We start with $X=P^2,$ its graph is the single vertex.
$$\circ$$
\vskip -0.6cm
$$-2$$

Blowing up a point, we get

$$\circ \!\! - \!\!\! - \!\!\! - \!\!\! - \!\! \circ$$
\vskip -0.6cm
$$-1\ \  -2$$

Blowing up another point on the pullback of $(\infty),$ we get 

$$\circ \!\! - \!\!\! - \!\!\! - \!\!\! - \!\! \circ \!\! - \!\!\! - \!\!\! - \!\!\! - \!\! \circ $$
\vskip -0.6cm
$$-1\ \  -2 \ \ -1$$

Blowing up a point on a newly blown up curve, we get 

$$\circ \!\! - \!\!\! - \!\!\! - \!\!\! - \!\! \circ \!\! - \!\!\! - \!\!\! - \!\!\! - \!\! \circ \!\! - \!\!\! - \!\!\! - \!\!\! - \!\! \circ  $$
\vskip -0.6cm
$$-1\ \  -2 \ \ -1 \ \ \ \  0$$

Then we blow up the intersection of the last two curves and get

$$\circ \!\! - \!\!\! - \!\!\! - \!\!\! - \!\! \circ \!\! - \!\!\! - \!\!\! - \!\!\! - \!\! \circ \!\! - \!\!\! - \!\!\! - \!\!\! - \!\! \circ  \!\! - \!\!\! - \!\!\! - \!\!\! - \!\! \circ  $$
\vskip -0.6cm
$$-1 \ \  -2 \ \ -1 \ \ \ -1 \ \ \ \ \  0\ $$

Blowing up another intersection point, we get 

$$\circ \!\! - \!\!\! - \!\!\! - \!\!\! - \!\! \circ \!\! - \!\!\! - \!\!\! - \!\!\! - \!\! \circ \!\! - \!\!\! - \!\!\! - \!\!\! - \!\! \circ  \!\! - \!\!\! - \!\!\! - \!\!\! - \!\! \circ  \!\! - \!\!\! - \!\!\! - \!\!\! - \!\! \circ  $$
\vskip -0.6cm
$$-1 \ \  -2 \ \ -1 \ \ \ -2 \ \ \  -1 \ \ \ \ \  0\ $$

\vskip -0.1cm
Blowing up another point, we get the following graph
\vskip -0.2cm

$$\circ \!\! - \!\!\! - \!\!\! - \!\!\! - \!\! \circ \!\! - \!\!\! - \!\!\! - \!\!\! - \!\! \circ \!\! - \!\!\! - \!\!\! - \!\!\! - \!\! \circ  \!\! - \!\!\! - \!\!\! - \!\!\! - \!\! \circ  \!\! - \!\!\! - \!\!\! - \!\!\! - \!\! \circ \!\! - \!\!\! - \!\!\! - \!\!\! - \!\! \circ  $$
\vskip -0.6cm
$$\ 0\ \ \ -1 \ \  -2 \ \ \ -1 \ \ \ -2 \ \ \  -1 \ \ \ \ \  0\ $$

Finally, blowing up four more points (in any order) we get the following:

\hskip 2.37cm $1\circ$ \hskip 6.15cm $\circ 1$

\vskip -0.2cm

\hskip 2.6cm $|$ \hskip 6.2cm $|$

\vskip -0.2cm

\hskip 2.6cm $|$ \hskip 6.2cm $|$

\vskip -0.7cm

$$\circ \!\! - \!\!\! - \!\!\! - \!\!\! - \!\! \circ \!\! - \!\!\! - \!\!\! - \!\!\! - \!\! \circ \!\! - \!\!\! - \!\!\! - \!\!\! - \!\! \circ  \!\! - \!\!\! - \!\!\! - \!\!\! - \!\! \circ  \!\! - \!\!\! - \!\!\! - \!\!\! - \!\! \circ \!\! - \!\!\! - \!\!\! - \!\!\! - \!\! \circ  $$
\vskip -0.7cm
$$0\ \ -1 \ \  -2 \ \ -1 \ \  -2 \ \  -1 \ \ \ \  0$$

\vskip -0.65cm

\hskip 2.6cm $|$ \hskip 6.2cm $|$

\vskip -0.2cm

\hskip 2.6cm $|$ \hskip 6.2cm $|$

\vskip -0.2cm

\hskip 2.37cm $1\circ$ \hskip 6.15cm $\circ 1$

For most of the exceptional curves, one can easily recover their self-intersection from the graph, using the adjunction formula:

$(K_{Z}+E_i)E_i=-2,$ so $\bar{K}_{Z}\cdot E_i=-2+\# (E_j \ {\textrm{adjacent to }}E_i)$

Thus, if $\bar{K}_{Z} = \sum a_i E_i,$ we have 
$$a_iE_i^2 + \sum \limits_{E_j \ adj. \ E_i} a_j = -2 + \# (E_j \ {\textrm{adjacent to }}E_i)$$

So if $a_i\neq 0,$ $E_i^2$ can be easily calculated.

However, when $a_i=0,$ it is not that easy. One can see in the above example, the left curve with $a_i=0$ has self-intersection $(-3),$ while the right one has self-intersection $(-4),$ despite the symmetry of the graph. One can remedy this situation by keeping track of the strict pullback of infinity. We do not need it in this paper, and the details are left to an interested reader.

Note that the subgraph of vertices with negative labels is connected. It is separated from the ``positive" vertices by the ``zero" vertices. Moreover, the ``zero" vertices are only connected to vertices with labels $(-1)$ or $1$.

Now we are going to make use of the map $\phi .$ The main idea is to use the adjunction formula for $\phi$ to get a formula for $\bar{K}_{Z}.$

Recall from the Introduction the four types of curves $E_i.$ For every curve of type $1$ or $3$ denote by $f_i$ the degree of the map onto its image and by $r_i$ the ramification index. Denote by $L$ the class of the line on $Y=P^2.$

\begin{Proposition} There exist integers $b_i$ for the curves $E_i$ of types $2$ and $4$ such that
$$\bar{K}_{Z} = \phi ^* (-2L) + \sum \limits_{type(E_i)=3} r_iE_i + \sum \limits_{type(E_i) =2 or 4} b_iE_i$$
\end{Proposition}

{\bf Proof.} Consider the differential $2$-form $\omega $ on $Y=P^2$ that has the pole of order $3$ at $(\infty)$ and no other poles or zeroes. Because $\phi $ is unramified on the $A^2 \subset X,$ there is a differential form on $Z,$ such that its divisor of zeroes and poles is $\phi ^*(-3L) + \sum _i c_iE_i,$ where $c_i$ can be calculated locally at a general point of each $E_i.$ 

Notice that for the curves $E_i$ of types $1$ and $3$, $c_i=r_i-1,$ and 
$$\phi ^* (L) = \sum \limits_{type (E_i)=1} r_iE_i + \sum \limits_{type (E_i)=2} e_iE_i$$
for some $e_i$. Thus,
$$\bar{K}_{Z} = K_{Z} + \sum E_i = \phi ^* (-3L) + \sum \limits_{type(E_i)=1 or 3} r_iE_i + \sum \limits_{type(E_i)=2 or 4} (c_i+1)E_i=$$
$$=\phi ^*(-2L) + \sum \limits_{type (E_i)=3} r_i E_i +\sum \limits_{type (E_i) =2 or 4} b_i E_i$$
 $\square$

Note that because $E_i$ are independent in the Picard group of $Z,$ the above representation of $\bar{K}_{Z}$ is unique and must match with the labeling of the graph of $E_i$. As a corollary, we have the following observation.

\begin{Proposition}
1) Any curve of type $1$ has a negative even label.

2) Any curve of type $3$ has a positive label.
\end{Proposition}

{\bf Proof.} Note that $\phi ^*(-2L)$ only involves curves of type $1$ and $2$.  $\square$
 
Additionally, the union of curves of type $1$ and $2$ must be connected, as a specialization (set-theoretically) of a pullback of a generic $L$ on $Y=P^2.$ This means that the corresponding subgraph is connected. 

Every curve of type $3$ must intersect with one of the curves of type $1$ or $2$, while the curves of type $4$ do not intersect with curves of type $1$ or $2$. (This follows from the projection formula of the intersection theory: if $E$ is  a curve on $Z,$ $E\cdot \phi^* (L) = (\phi_* E) \cdot L.$)

On the other hand, a type $3$ curve cannot intersect a type $1$ curve, because negative and positive labels are never adjacent. Because the graph of the exceptional curves on $Z$ is a tree, no two curves of type $3$ intersect with each other. Putting this all together, we must have the following. The tree of curves on $Z$ has a connected subtree containing all curves of type $1$ and $2$. Some of the vertices of this subtree may have one or more curves of type $3$ connected to them. Then some of these type $3$ curves may have trees of type $4$ curves connected to them. Additionally, no two curves of type $1$ are adjacent, and the subtree of curves of type $1$ and $2$ contains the connected subtree of curves with negative labels.

\begin{Proposition} $\pi ^{-1} (\infty )$ is of type $1$ or $2$.
\end{Proposition}

{\bf Proof.} One can prove it using the above description of the graph of exceptional curves, but there is also the following direct geometric argument. The pullbacks of lines on $X=P^2$ form a family of rational curves $C$ on $Z$ that intersect $\pi ^{-1} (\infty)$ at a generic point and do not intersect any other exceptional curves. Consider $\phi (C)$ for a generic $C.$ If $\pi ^{-1} (\infty) $ is of type $3$ or $4$ then $\phi (C) \subseteq A^2 \subset Y.$ The curve $C$ is proper and $A^2$ is affine, so $\phi (C)$ is a point, which is impossible.  $\square$

\vskip 0.3cm

Until now, the variety $Z$ was an arbitrary resolution at infinity of the original rational map. But we can put an additional restriction on it, to avoid unnecessary blow-ups.

\begin{Definition} If a curve is obtained by blowing up the intersection of two curves, we call these curves its parents. If a curve is obtained by blowing up a point on one of the curves, this curve is called its parent. The original line at infinity has no parents. Note that other curves may be created afterwards that separate the curve from one or both of its parents. 
\end{Definition}

\begin{Definition} For a given curve $E$, the set of its ancestors $A(E)$ is the smallest set $S$ of the exceptional curves that contains its parent(s) and  has the property that it contains the parents of every curve in $S$. Note that this set is empty if $E$ is the original line at infinity. Otherwise, it consists of the original line at infinity and all curves that have to be created before $E$. 
\end{Definition}

\begin{Definition}
A curve $E_i$ on $Z$ is called {\bf final} if there is a sequence of blow-ups from $X$ to $Z$ such that $E_i$ is blown up last. Equivalently, a curve is final if it is not a parent to any curve.
\end{Definition}

Note that there may be more than one final curve, and $\pi ^{-1} (\infty)$ is never final. In what follows, $E_i$ is one of the exceptional curves on $Z$.

\begin{Proposition}
Suppose that when $Z$ was created, $E_i$ was created after all of its neighbors in the graph (i.e. all adjacent vertices). Then $E_i$ is a final curve.
\end{Proposition}

{\bf Proof.} Instead of creating $E_i$ at its due time we can change the order of blow-ups and create it at the last step of the process, without changing anything else.  $\square$
  
\begin{Proposition}
Suppose $a_i=a(E_i)\geq 2$ and it is the largest  label among all its neighbors. Then $E_i$ is final.
\end{Proposition}

{\bf Proof.} We will prove that $E_i$ was created after all its neighbors. First of all, no neighbor of $E_i$ can be a blow-up of a point on $E_i,$ because its label would have been $a_i+1.$ If it were a blow-up of a point of intersection of $E_i$ and some $E_j,$ then $E_i$ and $E_j$ were adjacent before the blow-up. Negative curves are never adjacent to the positive curves and zero curves are only adjacent to curves with labels $1$ of $-1.$ Thus, $a_j \geq 1.$ So the label of the new curve is $a_i+a_j\geq a_i+1>a_i.$  $\square$

Note that no two curves with the same label $a_i\geq 2$ can be adjacent, by Proposition 3. So every local maximum $a_i\geq 2$ is a strict maximum. 

\begin{Proposition}
If $a_i=1,$ then $E_i$ is final if and only if it either has only one neighbor, with label $0$, or exactly two neighbors, with labels $1$ and $0.$ 
\end{Proposition}

{\bf Proof.} A curve with label $1$ can be created either by a blow-up of a point on a curve with label $0$ or by a blow-up of an intersection of a curve with label $0$ and a curve with label $1$. Once created, it will be final if and only if no other curve is blown up as its neighbor. The rest is easy and is left to the reader. $\square$ 

The above two propositions allow us to easily spot the final curves in the positive part of the graph of curves. Our interest in the final curves stems from the following. If one of the final curves on $Z$ is of type $2$ or $4$, then it can be contracted, using the $\phi -$relative MMP to get another $Z,$ with two maps to $X$ and $Y$ and a smaller Picard number.

\begin{Definition}
We call $Z$ {\bf minimal} if all of its final curves are of type $1$ or $3$.
\end{Definition} 

\begin{Proposition}
If a counterexample to the Jacobian Conjecture exists, it can be obtained using a minimal $Z.$
\end{Proposition}

{\bf Proof.} Take $Z$ with smallest possible Picard number. If it is not minimal, it can be created in such a way so that some curve of type $2$ or $4$ is blown-up last. Using MMP relative to $\phi,$ it can be blown down, maintaining the morphisms, and creating a counterexample to the Jacobian Conjecture with smaller Picard number.  $\square$

From now on, $Z$ will always be minimal.

\begin{Proposition} Suppose $E$ is a curve of type $3$ on $Z.$ Suppose $E_0$ is the curve of type $2$ it is adjacent to. Then the tree on the other side of $E$ is a line $E-E_1-...-E_k,$ where $E_1,...E_k$ are of type $4$.
\end{Proposition}

{\bf Proof.} The label of $E$ is positive.  All curves $E_1,...E_k$ ``on the other side" of $E$ are of type $4$. They must be ancestors of some curve of type $3$, so they are all ancestors of $E$. If the connected component of the graph obtained from $\Gamma$ by removing $E$ is not a line, there would have to be another  final curve curve there, which is impossible.  $\square$

\section{Other Varieties and Further Analysis}

We start with the theorem that shows that type $3$ curves must exist in a counterexample to the Jacobian Conjecture. Note that the type $3$ curves are called ``di-critical  components" in \cite{Domrina}, \cite{DomrinaOrevkov}, and this fact is well known and can be easily proven by a topological argument. So the main purpose of our proof is to show an easy application of our method before proceeding to the more intricate questions.

\begin{Theorem} Suppose $Z$ and $\phi$ provide a counterexample to the Jacobian Conjecture. Then $Z$ contains a curve of type $3$, where $\phi$ is ramified. 
\end{Theorem}

{\bf Proof.} 
Consider a generic line $L$ on the target variety $Y=P^2.$ The curve $C=\phi^{-1}(L)$ is smooth and irreducible (``Bertini's theorem"). Moreover, we can assume that for all but finitely many lines $L'$ that only intersect $L$ ``at infinity",  $C'=\phi^{-1}(L')$ is smooth and irreducible. We can also assume that $L$ does not pass through the images of the exceptional curves of types $2$ and $4,$ so $C$ does not intersect these curves on $Z.$ Suppose that the genus of $C$ is $g$, the map $H=\phi_{|_{C}}C\to L $ has degree $n$ and the number of points of $C$ ``at infinity" is $k$. (There is a special point $\infty$ on $L$, the only one not lying in $A^2.$ The number $k$ is the number of points of $C$ mapped to it, in a  set-theoretic sense.) Because the map $\phi$ is only ramified at the exceptional curves of $X,$ the map $H$ could only by ramified at these $k$ points at infinity. By Hurwitz formula, we have 
$$2g-2=-2n+r,$$
where $r$ is the total ramification at infinity. We have $g\geq 0,$ $n\geq 1$ and $r\leq n-k.$ So

$$-2\leq 2g-2 \leq -2n+n-k = -n-k \leq -2 $$
Thus all the inequalities above are equalities, $g=0,$ $n=1,$ and $k=1$. Because $n=1,$ the map $\phi$ is birational. For the birational maps the Jacobian Conjecture is well known (see, e.g. \cite{AbhyankarTataLectures}).  $\square$

Now we want to make further use of the morphism $\phi : Z \to Y.$ We decompose it into a composition of two morphisms, birational and finite (Stein factorization):
$$Z \longrightarrow W \longrightarrow Y$$ 

Here the first morphism is birational and denoted $\tau$, and the second one is finite and denoted $\rho.$ 

The surface $W$ is  normal. In what follows, we use the intersection theory for normal surfaces due to Mumford. Suppose $K_{W}$ is its canonical class, as the Weil divisor class modulo numerical equivalence. We define the augmented canonical class $\bar{K}_{W}=K_{W}+\sum E_i$.

\begin{Proposition} In the situation and notation described above,
 $$\bar{K}_{W}= \rho^*(-2L) +\sum \limits_{type(E_i)=3} r_iE_i.$$
\end{Proposition}

{\bf Proof.} The curves $E_i$ on $W$ are exactly the images of curves of types 1 and 3 on $Z.$ By adjunction, we have:
$$K_{W} = \rho^* K_Y + \sum (r_i-1) E_i,$$
where $r_i$ is the ramification index, and $E_i$ are the images of the curves $E_i$ of types $1$ and $3$ on $Z.$
$$\bar{K}_{W}=\rho^*(-3L)+\sum r_iE_i = \rho^*(-2L) +\sum \limits_{type(E_i)=3} r_iE_i.$$
 $\square$

\vskip 0.5cm

Note that if one denotes by $\bar{K}_Y$ the class of $K_Y + (\infty ),$  then $\rho^*(-2L)$ in the above Proposition is $\rho^*(\bar{K}_Y).$
The next theorem is very important. It will be further strengthened at the end of this section.

\begin{Theorem} (Big Ramification Theorem)

Suppose $Z$ is a counterexample to the Jacobian Conjecture. Then on $W$ the ``di-critical log-ramification divisor"
$$\bar{R}=\sum\limits_{E_i\subset W, type (E_i) =3} r_i  E_i $$
intersects positively with all exceptional  curves of type $3$. As a corollary,
 $\bar{R}^2 >0.$
\end{Theorem}

{\bf Proof.} Suppose  that $C =E_i$ is a curve of type $3$ on $W.$ Suppose $d_i$ is the degree of $\rho(C)$. Suppose $\tau : Y_2 \to W$ is a minimal resolution of singularities of  $W$. Then
$$\bar{K}_{W}E_i=(K_{W}+E_i)E_i + \sum \limits_{j\neq i} E_iE_j,$$
where $E_j$ are curves of type $1$ or $3$. Note that because $E_i$ intersects at least one curve of type $1$, $\bar{K}_{W}E_i>(K_{W}+E_i)E_i.$
 Lifting up to $Y_2,$ we get
$$(K_{W}+E_i)E_i =(\tau^*(K_{W})+\tau^*(E_i))\tau^{-1}(E_i) \geq (K_{Y_2}+\tau^{-1}(E_i))\tau^{-1}(E_i) \geq -2. $$
So for all $i$ $\bar{K}_{W}E_i>-2.$

Therefore,
$$\bar{R}\cdot E_i = 2\phi^*(L)\cdot  E_i+ \bar{K}_{W}\cdot E_i   >2f_id_i-2 \geq 0$$
 $\square$ 

\begin{Corollary} The curve $\pi ^{-1} (\infty)$ is of type $2.$
\end{Corollary}

{\bf Proof.} By Proposition 2.6, it is of type $1$ or $2$. If it is of type $1$, then it is not included in $ \sum \limits_{type(E_i)=3}  r_i\tau^*E_i$. So $ \sum \limits_{type(E_i)=3} r_i\tau^*E_i$  consists of curves contractible by $\pi .$ But 
$$(\sum \limits_{type(E_i)=3} r_i \tau^*E_i)^2 = (\sum \limits_{type(E_i)=3} r_iE_i)^2 >0, $$
contradiction.  $\square$

Note that every curve of type $3$ on $Z$ intersects the union of curves of type $2$ at exactly one point, and does not intersect curves of type $1.$ When the curves of type $2$ are contracted, on $W,$ every curve of type $3$ intersects the union of curves of type $1$ at exactly one point.

\begin{Proposition}
For every curve $E_i$ of type $3$ on $W$ the point above is $\tau(\pi^{-1} (\infty))$.
\end{Proposition} 

{\bf Proof.} Suppose there is a point $w\in W$ on the union of type $1$ curves, which is not $\tau (\pi^* (\infty))$ and which has some type $3$ curves passing through it. Define
$$\bar{R}_w= \sum  \limits_{w\in E_i, type(E_i)=3} r_iE_i$$
Because the curves of type $3$ not passing through $y$ cannot intersect any components of $\bar{R}_w,$ we have $\bar{R}_w^2 = \bar{R}_w \cdot \bar{R}$. By Theorem 3.2, this implies that $\bar{R}_w^2>0.$ Like in the proof of Corollary 3.1, $\tau^*(\bar{R}_w)$ consists of curves contractible by $\pi,$ which is impossible. $\square$

\begin{Proposition}
On $W,$ all exceptional curves contain $\tau (\pi ^{-1} (\infty))$ and there are no other points of intersection.
\end{Proposition}

{\bf Proof.} By the proposition above, every curve of type $3$ contains $\tau (\pi ^{-1}(\infty))$ and this is its only point of intersection with other exceptional curves. Now consider a curve $E_i$ of type $1$. Suppose it does not contain $\tau (\pi^{-1}(\infty))$. Then it does not intersect any of the curves of type $3$ on $W.$ 

On $W$ we have:
$$\bar{K}_{W} \cdot E_i \geq  (K_{W}+E_i) E_i \geq -2 $$

On the other hand,
$$\bar{K}_{W} \cdot E_i = (-2\rho^*(L)+\bar{R} ) \cdot E_i = -2\rho^*(L) \cdot E_i \leq -2 $$

The inequalities above become equalities only if $E_i$ intersects no other curves and is smooth. This would make it the only curve of type $1$ on $W$, which would then have to intersect with some curves of type $3,$ contradiction.  $\square$

\vskip 0.5cm

Thus, we know that every curve of type $1$ on $W$ contains $\tau (\pi ^{-1} (\infty)).$ We now look at the graph of curves on $Z.$ The curves of type $2$ that are mapped to $\tau (\pi ^{-1} (\infty))$ form a connected subgraph, containing $\pi ^{-1} (\infty).$ Every curve of type $1$ or $3$ is attached to this subgraph. On ``the other side" of each curve of type $3$ there may be a single chain of curves of type $4$, and on ``the other side" of each curve of type $1$ there may be a single chain of curves of type $2$. Note that all of these "other side" curves must be created before the corresponding type $3$ or type $1$ curves. When mapped to $W,$ the curves of type $1$ and $3$ intersect at $\tau (\pi ^{-1} (\infty))$ and nowhere else.

One can restrict the structure of the possible counterexamples even further.
\begin{Theorem} In any counterexample to the Jacobian Conjecture there are no curves on ``the other side" of the curves of type $1$.
\end{Theorem}

{\bf Proof.} Consider a curve of type $1$, $E,$ on $Z$. Suppose the ramification index at $E$ is $r$. Then the coefficient of $\phi ^* (L)$ in $E$ is $r$, and the coefficient of $\bar{K_{Z}}$ is $(-2r).$ Consider the divisor class $D=\bar{K_{Z}}+2\phi* (L)= ...+0\cdot E +x_1E_1+...+x_kE_k,$ where $E_1,...E_k$ are the curves on $Z$ ``on the other side" of $E$.
We know that $D$ intersects by zero with $E_1,...,E_{k-1}$. It intersects by $-1$ with $E_k$. We formally add another vertex to the graph, $"E_{k+1}"$ and set the coefficient of $D$ at it to be $1$. (Note that we are not blowing up any points and $E_{k+1}$ does not have any geometric meaning). We now have a chain $E, E_1, ..., E_{k}, E_{k+1}$ and a divisor $D'= 0\cdot  E +x_1E_1+...+x_kE_k +1\cdot E_{k+1},$ such that $D'$ intersects by zero with all $E_1, E_2,...,E_k.$ Because the self-intersections of all $E_i,$ $1\leq i \leq k,$ are less than or equal to $-2,$ the coefficients $x_i$ must form a concave up chain between $0$ and $1$, contradicting their integrality. (Here is a more formal argument. Suppose at least one of the $x_i,$ $1\leq i \leq k,$ is not positive. Then consider the minimum of $x_i$, obtained at $x_j$, such that $x_{j+1}> x_j$, where formally $x_0=0, x_{k+1}=1$. Then $D'\cdot E_j \geq x_{j-1}+x_{j+1} +2x_j >0 ,$ contradiction.  Suppose the maximum of $x_i,$ $1\leq i \leq k,$ is greater than or equal $1$ and is obtained at $x_j,$ where $x_{j-1}<x_j.$  Then $D'\cdot E_j \leq x_{j-1}+x_{j+1} - 2 x_j <0,$ contradiction. Thus all $x_i$ are strictly between $0$ and $1$, which is impossible because they are integers.)  $\square$

As a corollary of these observations, we get a rather detailed description of the structure of $W.$
\begin{Theorem} Suppose $W$ is defined as above for a counterexample to the Jacobian Conjecture, $E_i$ are images of the exceptional curves on it. Then  $W\setminus \cup_i E_i$ is isomorphic to the affine plane. There exist distinct points $A=\tau( \pi^{-1}(\infty))$ and $A_i \in E_i$ (at most one for each $E_i$ of type $3$) so that $W$ is smooth outside of them, $A$ is a normal singularity, $A_i$ are cyclic quotient singularities, all the curves $E_i$ pass through $A,$ do not intersect elsewhere and are smooth in the nonsingular part of $W$. For all exceptional curves $E_i$ that do not contain $A_i,$ $E_i\setminus \{A\}$ is isomorphic to the affine line. For al curves $E_i$ that contain $A_i,$   $E_i\setminus \{A,A_i\}$ is isomorphic to the algebraic torus (affine line with a removed point).
\end{Theorem}

{\bf Proof.} Most of the statements have already been proven. To finish the proof, note the following. The map $\tau$ from $Z$ to $W$ contracts all curves of types $2$ and $4$, and no curves of types $1$ and $3$. The curves of type $2$ form a subtree on $Z,$ so they are contracted to one singular point $A$. Note that all curves of type $1$ and $3$ on $Z$ intersect this subtree at exactly one point. Some curves of type $3$ have one chain of curves of type $4$ attached to them, that get contracted into a cyclic quotient singularity.  $\square$

Note that the above theorem restricts greatly the restriction of the map $\phi$ to the exceptional curves of types $1$ and $3$. 

\begin{Theorem} 1) For all curves $E_i$ of type $1$, $f_i=1.$

2) For all curves $E_i$ of type $3$ either $f_i=1$ or the restriction to $E_i$ of the map from $Z$ to $Y$  is isomorphic to the composition of a map $(x\mapsto x^{f_i} ):P^1\to P^1$ and a generically one-to-one  map from $P^1$ to a possibly singular rational curve (the normalization map for $\phi(E_i)$).
\end{Theorem}

{\bf Proof.} This follows from the fact that the restriction  to $E_i$ of the map from $Z$ to $Y$ can only be ramified at the points of intersections of $E_i$ and other exceptional curves, and the classification of self-maps of  the projective line that are ramified at one or two points.  $\square$

For the following theorem, we need to introduce additional notation.

\begin{Definition}
Suppose $Z,W,Y$ and $\phi :W\to Y$  are as above. For a type $3$ curve $E_i,$ denote by $F_i$ its  image on $Y$, as a reduced irreducible divisor (a possibly singular rational curve). Then 
$$\rho^*(F_i)=r_iE_i+G_i,$$
where $G_i$ is an effective Weil divisor. Its irreducible components are curves in $A^2$ that are mapped to $F_i.$ We will call these curves coexceptional.
\end{Definition}

\begin{Theorem} Suppose $E_i\subset W$ is a type $3$ curve that contains a cyclic quotient singularity $A_i.$ Then some coexceptional curve from $G_i$ contains $A_i$.
\end{Theorem}

{\bf Proof.} The proof is very similar to the proof of Theorem 3.3. Suppose the support of $G_i$ does not contain $A_i.$ To simplify the notation, denote $E_i$ by $E$; suppose $E_1, E_2,...,E_k$ are the curves of type $4$ that are mapped to $A_i,$ with $E_1$ intersecting $E$. Suppose $\phi :Z \to Y$ is our map. Consider on $Z$  the divisor $D=\bar{K}_{Z}-\phi^*(F_i).$  We can write $D$ as a linear combination of exceptional curves and the strict pullbacks of the coexceptional curves. Because the coefficient of $E$ in this linear combination is zero, and the only curves that contribute to the intersection of $D$ with $E_1, ..., E_k$ are $E_1,...E_k,$ we get a linear combination with integer coefficients $x_1E_1+...+x_kE_k$ that intersects by zero with $E_1,...,E_{k-1}$ and by $(-1)$ with $E_k.$ We can now follow verbatim the argument in Theorem 3.3 to get a contradiction.  $\square$

We end this section with a strengthening of Theorem 3.2.

\begin{Theorem} (Ample Ramification Theorem)

Suppose $Z$ is a counterexample to the Jacobian Conjecture. Then on $W$ the ``di-critical log-ramification divisor"
$$\bar{R}=\sum\limits_{E_i\subset W, type (E_i) =3} r_i  E_i $$ is ample.
\end{Theorem}

{\bf Proof.} The surface $W$ is rational and, therefore, $\Q$-factorial. So the Weil divisor $\bar{R}$ is $\Q$-Cartier. It is effective and  by Theorem 3.2 it intersects positively with all of its irreducible components. So by the Nakai-Moishezon criterion it is enough to show that it intersects positively with all irreducible curves $C$ on $W$ that are not the exceptional curves of type $3$.

If $C$ is a curve of type $1$ on $W$, then it intersects the support of $\bar{R}$ at $A=\phi  (\pi^{-1} (\infty)),$ so $C\cdot \bar{R}  >0.$ 

If $C$ is any other curve on $W$, that does not intersect positively with $\bar{R}$, it does not intersect with any curves of type $3$ and it must intersect with at least one curve of type $1$ (because it cannot be contained entirely in the affine plane). By the Hodge Index Theorem and Theorem 3.2, $C^2 <0$. Note also that $C$ does not pass through $A,$ so it intersects with at least one curve of type $1$ at a smooth point. Because the $\bar{K}$-labels of all curves of type $1$ are at most  $-2,$ $K_{W}\cdot C \leq -3.$ Therefore, $(K_{W}+C)\cdot C \leq -3 < -2,$ which is impossible.  $\square$

\section{Adjunction Formulas and Inequalities}

As the above discussions clearly indicate, one of the basic objects that we need to consider are algebraic surfaces with a collection of reduced irreducible curves on them. The natural framework for this is provided by the notion of varieties with boundaries. They have been studied extensively as part of the Minimal Model Program, and most of our results in this section are not new, and hold in a much more general setting. A good introduction to the more general theory can be found in \cite{Kollar_pairs}. We start by recalling several standard definitions, but this section is directed at the readers with some background in modern birational geometry.  All varieties are over $\C$, all surfaces are normal unless otherwise specified and $\Q$-factorial whenever needed. Please note that the notation in the beginning of this section is independent of the notation in the rest of the paper.

\begin{Definition} A log-surface (or surface with boundary) is a pair $(X,D)$ where $X$ is an algebraic surface, and $D=\sum_i D_i$ is a sum of prime divisors (i.e. reduced irreducible curves) on $X.$ The augmented canonical class of $(X,D)$ is  $\bar{K}=K+D,$ where $K$ is the canonical class of $X.$ The divisor $D$ is called the boundary, and $D_i$ are the components of the boundary. We will often omit $D$ from the notation, when it is clear from the context what it is (usually, a complement of $A^2$ in $X$).
\end{Definition}

\begin{Definition} A resolution of singularities of a log-surface $(X,D)$ is a smooth algebraic surface $Y$ with a birational map $\tau :Y \to X$ such that the exceptional divisors (curves) of $\tau$ and the components of $\tau^{-1}(D)$ have simple normal crossings (i.e. all the irreducible components are smooth and intersect transversally).
\end{Definition}

Note that whenever we have a birational map $\tau:Y \to X,$ we take for the boundary on $Y$ the union of the components of $\tau^{-1}(D)$ and the exceptional divisors.  A resolution as above always exists, moreover it can be obtained by a sequence of blow-ups of the surface or the irreducible components of the boundary or the points where the components of the boundary do not intersect transversally.

\begin{Definition} For a birational map $\tau :Y \to X$, the adjunction formula is an equality of $\Q-$Cartier divisors
$$\bar{K}_Y = \tau^*\bar{K}_X +\sum \limits_i a_iE_i,$$
where the sum is over the exceptional divisors of $\tau$ and $a_i$ are rational numbers called log-discrepancies. 
\end{Definition}

\begin{Definition} A log-surface (X,D) has log-canonical singularities iff for one (any) resolution of singularities $Y$ as above all log-discrepancies are nonnegative.
\end{Definition}

\begin{Example} If $X$ is smooth and $D$ has simple normal crossings, then $(X,D)$ has log-canonical singularities. 
\end{Example}

\begin{Example} Suppose a point $P\in X$ is a normal surface singularity and $D_1$ and $D_2$ are curves on $X$ containing $P$. Suppose the dual graph of exceptional curves on the minimal resolution $\pi : \tilde{X} \to X$ of $X$ consists of a chain of rational curves  $E_1,...E_k$  (i.e. all $E_i$ are isomorphic to $P^1$ and $E_i$ intersects $E_{i+1}$ transversally for $1\leq i\leq k-1$ with no other points of intersection). Suppose that the strict pullback of $D_1$ intersects $E_1$ transversally and the strict pullback of $D_2$ intersects $E_k$ transversally, so that $\pi^{-1}D_1,E_1,...E_k,\pi^{-1}D_2 $ is a chain of rational curves.
 Denote $D=D_1+D_2$. Then $(X,D)$ is log-canonical in the neighborhood of $P$. (Note that in this case the singularity of $X$ at $P$ is locally-analytically a cyclic quotient). Note that in this case the resolution $\tau : Y\to X$ is log-crepant (i.e. all log-discrepancies are zero).
\end{Example}

The following theorem is the result of the relative Minimal Model Program. 

\begin{Theorem} Suppose $(S,D)$ is a log-surface, and $\tau : V \to S$ is its resolution of singularities. Then the map $\tau$ can be decomposed into a composition of two maps, $\chi:V\to S_1$ and $\mu: S_1 \to S$, such that

1) The log-surface $(S_1,D_1)$ has log-canonical singularities;

2) The augmented canonical class $\bar{K}_{S_1}$ is $\mu $-ample (i.e. $\bar{K}_{S_1} E_i >0$ for all $E_i$ that are contracted to a point by $\mu$).

Here $D_1$ is, as usual, the union of components of  $\mu^{-1} D $ and the exceptional divisors of $\mu$. 
This $S_1$ with the morphism $\mu$ is called the log-canonical modification of $S$.
\end{Theorem}

{\bf Proof.} This can be done by contracting one-by-one exceptional curves of $\tau$ that have negative intersection with the augmented canonical class. See, e.g. the paper of Fujino \cite{Fujino}, section 3 for this theory in a much more general setting. $\square$

Because of the above theorem, it is important to classify all log-canonical singularities with boundary coefficients $1$. The general classification of log-canonical surface singularities is well known, see the work of Alexeev \cite{Alexeev_logcanonical}. We are most interested in the surface being log-canonical along the boundary $D,$ especially at the points of intersections of components of $D.$ 

\begin{Theorem} Suppose $P\in X$ is a point of intersection of two components of $D$.  Then $P$ is a log-canonical singularity if and only if it is of the kind described in the Example 4.1 or 4.2.
\end{Theorem}

{\bf Proof.} This follows from the classification of \cite{Alexeev_logcanonical}, and can also be obtained independently. $\square$

We now proceed to applications, and need to introduce a new definition.

\begin{Definition} Suppose $(S,D)$ is log-surface. Suppose $D_i$ is an irreducible component of $D.$ Then the valency of $D_i$ is the number of points in the normalization of $D_i$ that are mapped by the normalization map to the points of intersection of $D_i$ and some other $D_j$. We denote the valency of $D_i$ by $val_D(D_i)$ or simply $val(D_i)$ if $D$ is clear from the context.
\end{Definition}

The notion of valency generalizes to the singular surfaces the notion of the number of neighbors in the dual graph of exceptional curves. The following theorem  follows immediately from the adjunction formula for the curve $D_i$ on the surface $S$.

\begin{Theorem} Suppose $S$ is a smooth surface, and the boundary $D$ is a union of not necessarily smooth curves intersecting transversally at smooth points (i.e. $D$ has normal crossings but not necessarily simple normal crossings). Suppose $D_i$ is one of the irreducible components of $D$. As usual, $\bar{K}_S$ stands for $K_s+D.$ Then
$$\bar{K}_S \cdot D_i \geq  (2g-2) +val(D_i),  $$
and the equality occurs if and only if $D_i$ is smooth.
\end{Theorem}

The following theorem is an appropriate generalization of the above to arbitrary surfaces with boundaries. 

\begin{Theorem} (Adjunction Inequality) Suppose $(S,D)$ is a log-surface, $D_i$ is an irreducible component of the boundary. Suppose the genus of the normalization of $D_i$ is $g$. Then
$$\bar{K}_S \cdot D_i \geq (2g-2) +val(D_i)$$ 
Moreover the equality occurs if and only if all of the following conditions are satisfied:

1)  $(S,D)$ is log-canonical in the neighborhood of $D_i$;

2) $S$ and $D_i$ are smooth at all points $P$ on $D_i$ that do not lie on other boundary components.
\end{Theorem}

{\bf Proof.} Suppose $(S',D')$ is the log-canonical modification of $(S,D),$ with the morphism $\mu :S' \to S.$ Because $\bar{K}_{S'}$ is $\mu-$ample, all log-discrepancies $a_i$ in the adjunction formula
$$\bar{K}_{S'} = \mu ^* \bar{K_S} + \sum_i a_iE_i$$
are positive. (Note that if $(S,D)$ is log-canonical, then $S'=S$ and then there are no exceptional curves $E_i$).

Denote $D'_{i} = \mu^{-1} D_i$.  Note that if $S'$ is not equal to $S$ in the neighborhood of $D_i$, then
$$\bar{K}_{S'} \cdot D'_{i} < \bar{K}_{S} \cdot D_i $$
 It is also easy to see that $val(D'_{i})=val(D_i)$. So to prove the desired inequality for $D_i$ it is enough to prove it for $D_{i}$.

Consider the minimal resolution of singularities $ \chi: V\to S'$. Above every point of intersection of $D_i$ with other boundary components we get a chain of rational curves with the end curves intersecting transversally with $D_i$. Define the boundary on $V$ to be the union of $\chi^{-1} D'$ and these exceptional curves above the points of intersection of $D'_i$ and other components of $D'.$  Consider the adjunction formula:

$$\bar{K}_{V} = \chi^* \bar{K}_{S'} + \sum a_j E_j,$$
where $E_j$ are the exceptional curves of $\chi.$

Note that for the curves $E_j$ that are mapped to the points of intersection of $D'_i $ and other components of $D'$ the log-discrepancies $a_j$ are equal to zero. On the other hand, for all curves $E_j$ that are mapped elsewhere on $D'_i$ the discrepancies $a_i$ are strictly negative, because $ \chi^*K_{S'}- K_V$ is effective and $\chi^*D' -\chi^{-1} D' $ is a sum of $E_j$ with strictly positive coefficients. Multiplying by $\chi ^{-1}D'_i,$ we get
$$\bar{K}_V\cdot  \chi ^{-1}D'_i \leq  \bar{K}_{S'} \cdot D'_{i}   $$

Applying Theorem 4.3 to $V$ and $\chi^{-1}D'_i$, we get the desired inequality. Also, in order to have equality, several conditions need to be satisfied.

1) The surface $S'$ must equal $S$ in the neighborhood of $D'_i$, that is $(S, D)$ is log-canonical in the neighborhood of $D_i$. 

2) No exceptional curves are mapped by $\chi $ to a point on $D'_i$ that does not lie on other components of the boundary. 

3) The curve $\chi ^{-1}D'$ must be smooth. By the previous conditions this implies that $D$ is smooth outside of intersections with other boundary components.

Finally, it is easy to check that when these conditions are satisfied, all inequalities become equalities.
 $\square$

\begin{Definition} If the boundary $D$ on $S$ satisfies the conditions of the above theorem for all irreducible components, we say that $D$ has generalized simple normal crossings.
\end{Definition}

\begin{Remark} The number $(2g-2)+val(D_i)$ is the degree of the augmented canonical divisor on the normalization of $D_i,$ with the boundary being the pullback of the union of the other components of $D$ by the normalization map. 
\end{Remark}

We now return our attention to a hypothetical counterexample to the Jacobian Conjecture. Generalizing the construction of previous sections,  suppose $Y$ is a compactification of the target $A^2,$ $X$ a compactification of the origin $A^2;$ $\pi: Z \to X$, $\phi:Z\to Y$ is a minimal resolution of the rational map from $X$ to $Y$ and $Z\to W \to Y$ is the Stein factorization of $\phi$, with the first map denoted by $\tau$ and the second one by $\rho$. Because $\phi $ is unramified on $\pi^{-1} (A^2)$, all the ramification of $\rho$ occurs at the images by $\tau$ of the complement of $A^2$ on $Z$ (i.e. the exceptional curves). All these surfaces will be equipped with the boundaries that are complements of $A^2$ that they naturally contain. The irreducible components of these boundaries will be called exceptional curves. Note that all exceptional curves are rational.

In this general setting we can still distinguish four classes of exceptional curves on $Z$. However, it is important to keep track for curves of type $1$ what exceptional curve they are mapped to, and for curves of type $2$ what exceptional curve(s) contain their image. Of most interest to us are the curves of type $3$, i.e. the di-critical curves. Unlike the curves of type $1$ and $2$, they are essentially independent of the choices of the compactification. Just as before, the exceptional curves of type $2$ and $4$ get contracted by $\tau,$ so on $W$  exceptional curves of types $1$ and $3$ remain. For each curve $E_i$ of type $1$ or $3$ on $W$ we will denote by $F_i$ its image by $\mu$. We denote by $r_i$ the corresponding ramification index and by $f_i$ the degree of the restriction of $\mu$ to $E_i.$ We will often denote the di-critical components by $R_i$ to distinguish them from other exceptional curves. We will call $\bar{R} = \sum r_i R_i$ the di-critical log-ramification divisor. Its importance lies in the fact that it connects the augmented canonical classes of $W$ and $Y$ as follows.

\begin{Theorem} (Adjunction Formula for $\rho$) Suppose $W,$ $Y$ and $\rho$ are as above. Then
$$\bar{K}_W= \rho^* \bar{K}_Y + \bar{R}$$
\end{Theorem}

{\bf Proof.} The argument is the same as in Proposition 3.1: because our polynomials map sends $A^2$ to $A^2,$ the preimage of every boundary component on $Y$ is the union of the boundary components on $W$. The only ramification not accounted for in this way is in the di-critical components; since the boundary on $W$ includes them, we get $\bar{R}$ on the right hand side. The equality should be understood as the rational equivalence of Weil divisors. $\square$

In many cases, the surface $Y$ will be smooth, with simple normal crossing boundary, but the argument works equally well when $Y$ is a log-canonical surface with generalized simple normal crossing boundary. It is in this more general setting that we state the following theorem.

\begin{Theorem} Under the above conditions, suppose $(Y, Y\setminus A^2)$ is a surface with generalized simple normal crossing boundary. Suppose $E_i$ is an exceptional curve on $W$ such that $E_i$ does not intersect any di-critical curves on $W.$ Then $(W, W\setminus A^2)$ has generalized simple normal  crossings in the neighborhood of $E_i$.
\end{Theorem}

{\bf Proof.} It is a standard matter to prove that $W$ is log-canonical in the neighborhood of $E_i$. However we need a slightly stronger result, and will use Theorem 4.4. for it. Suppose $\rho (E_i) = F_i$, the ramification index of $\rho$ at $E_i$ is $r_i$ and the degree of the restriction of $\rho$ to $E_i$ is $f_i.$ Multiplying the adjunction formula for $\rho$ by $E_i,$ we get
$$\bar{K}_W \cdot E_i = f_i \cdot (\bar{K}_Y \cdot F_i)$$
By Theorem 4.4 and Remark 4.1, $\bar{K}_Y \cdot F_i$ is the degree of the augmented canonical divisor on the normalization of $F_i$. The adjunction formula for the map from the normalization of $E_i$ to the normalization of $F_i$ (the Riemann-Hurwitz formula) implies that
$\bar{K}_W \cdot E_i$ is the degree of the augmented divisor on the normalization of $E_i$. (Note that the boundary there is exactly the pullback of the boundary on the normalization of $F_i$, and all ramification must occur there). Theorem 4.4. now implies the desired result. $\square$

We finish this section with some simple observations that describe explicitly the situation when the valency of $F_i$ or $E_i$ is small.

\begin{Theorem} Under the conditions of the above theorem, the following is true. 

a) If $val (F_i)=1,$ then $val(E_i)=1$ and $f_i=1.$

b) If $val (F_i)=2,$ then $val(E_i)=2$. If $f_i>1,$ then the map from $E_i$ to $F_i$ is ramified at two points.
\end{Theorem}

{\bf Proof.} As in the above theorem, 
$ \bar{K}_W \cdot E_i = f_i \cdot (\bar{K}_Y \cdot F_i)$.

Note that $E_i$ and $F_i$ are rational, so the above theorem and Theorem 4.4. imply that
$val(E_i)-2=f_i\cdot (val(F_i)-2)$. This implies the result, considering that the map from $E_i$ to $F_i$ can only be ramified at the points of intersection with the other boundary components. $\square$

\begin{Theorem} Under the conditions of Theorem 4.6, the following is true. 

a) If $val (E_i)=1,$ then $val(F_i)=1$ and $f_i=1.$

b) If $val (E_i)=2,$ then $val(F_i)=2$. If $f_i>1,$ then the map from $E_i$ to $F_i$ is ramified at two points.

c) If $val (E_i)=3,$ then $val(F_i)=3$ and $f_i=1.$
\end{Theorem}

{\bf Proof.} As in the proof of the previous theorem, $val(E_i)-2=f_i\cdot (val(F_i)-2).$ The results follow immediately. $\square$

\section{Determinants of Weighted Forests and Some Applications}

In this section we will construct and use a more sophisticated labeling on the graph of exceptional curves, using the determinants of the matrices related to the Gram matrix of the intersection form. This approach is related to the work of Domrina and Orevkov (cf. \cite{Domrina}, \cite{DomrinaOrevkov}) and some of our results can be derived from theirs and/or the work of Walter Neumann (\cite{Neumann}). However our approach is slightly different: we do not use the splice diagrams and we derive applications by combining this labeling with the $\bar{K}-$ labeling described above. 

In what follows, a weighted tree means a connected graph with no cycles, with weights attached to the vertices. A weighted forest is a graph with no cycles, with weights attached to the vertices. To each such graph $\Gamma$ and any ordering of the vertices we will associate a matrix $Q(\Gamma)$ as follows. 

$$Q(\Gamma)_{i,j}=\left\{
\begin{array}{l}
-a_{i},  \ if\ i=j, \\
-1, \  if\  i^{th}\  and\  j^{th}\  vertices\  are\  connected\  by\  an\  edge\\
0\  {\rm otherwise}
\end{array}
\right.$$

In applications, the forest will be obtained from a tree of exceptional curves on $Z$ by removing some vertices and edges, with weights being the self-intersection numbers. So the matrix $Q$ will be related to a Gram matrix of the minus-intersection form. Note that our choice of a matrix differs from the usual matrix of a weighted graph. This is done because the minus-intersection form is ``almost" positive-definite.

Suppose $\Gamma$ is a weighted forest. Denote by $E(\Gamma)$ and $V(\Gamma)$ the sets of edges and vertices of $\Gamma$ respectively. Denote by $d(\Gamma)$ the determinant of $Q(\Gamma)$ (note that it does not depend on the ordering of the vertices). For every subset $S$ of $E(\Gamma)$ denote by $Supp(S)$ the set of all vertices that are endpoints of some edge from $S$. A subset of $S$ will be called disjoint if no two edges share a vertex (i.e. $|Supp(S)|=2|S|$). Denote by $A(\Gamma)$ the set of all disjoint subsets of $E(\Gamma)$. The following lemma provides a useful formula for the determinant $d(\Gamma)$.

\begin{Lemma} In the above notations, if the weight of $v\in V(\Gamma)$ is $a_v$,
$$d(\Gamma) = \sum \limits_{S\in A(\Gamma)} \Big ( (-1)^{|S|} \cdot \hskip -0.3cm \prod \limits_{v\in V(\Gamma) \setminus Supp(S)} (-a_v) \Big)$$
\end{Lemma}
{\bf Proof.} Using ordering of $V(\Gamma),$ we identify it with $\{1,2,\dots, |V(\Gamma)| \}$. The determinant $d(\Gamma)$ is the sum of $|V(\Gamma)|!$ terms corresponding to permutations of $V(\Gamma)$. Note that such term is zero unless the corresponding permutation sends each $i$ to itself or to a vertex $j$ which is connected to $i$ by an edge. Because $\Gamma$ is a forest, the cycle decomposition of such permutation consists of fixed points and transpositions $(i,j)$, where $(i,j)$ are edges of the graph. Obviously, these edges must be disjoint, so the set of such permutations is in one-to one correspondence with $A(\Gamma)$. The formula follows, considering that the sign of such permutation is determined by the parity of the number of transpositions.  $\square$

We will use the following notations and definitions.

\begin{Definition}
1) For $P\in V(\Gamma)$ we denote by $\Gamma_P$ the weighted forest obtained from $\Gamma$ by removing $P$ and all edges involving $P$.  We define the determinant label $d_P$ of the vertex $P$ as $d(\Gamma_P).$

2) For $(P,Q)\in E(\Gamma)$ we denote by $\Gamma_{PQ}$ the weighted forest obtained from $\Gamma$  by removing  the edge $(P,Q).$   We define the determinant label $d_{PQ}$ of the edge $(P,Q)$ as $d(\Gamma_{PQ}).$

3) For a subset $\{P,...,Q\}$ of $V(\Gamma)$ we denote by $\Gamma_{\{P,...,Q\}}$ the weighted forest obtained from $\Gamma$ by removing all vertices in this subset and all edges involving them. We denote by $d_{P,...,Q}$ the determinant of $\Gamma_{\{P,...,Q\}}$. {\bf Note that for an edge $(P,Q)$ the number $d_{P,Q}$ should not be confused with $d_{PQ}$}. 

\end{Definition}

The following four properties of the determinants will be used extensively later. We will call them the multiplicativity property, the expansion by vertex formula,  the expansion by edge formula, and the weight increment formula.

\begin{Lemma} Suppose $\Gamma$ is a disjoint union of $\Gamma_1$ and $\Gamma_2$. Then $$d(\Gamma)=d(\Gamma_1)\cdot d(\Gamma_2)$$
\end{Lemma}

{\bf Proof.} With the appropriate ordering of the vertices, the matrix $Q(\Gamma)$ is block-diagonal, with blocks being $Q(\Gamma_1)$ and $Q(\Gamma_2)$. $\square$

\begin{Lemma} Suppose $P\in V(\Gamma)$.  Then
$$d(\Gamma)=-a_P\cdot d_P  - \sum \limits_{(P,Q) \in E(\Gamma)} d_{P,Q}$$
\end{Lemma}

{\bf Proof.} This follows from Lemma 4.1 by breaking up the terms in the sum for $d(\Gamma)$ into classes based on the way the vertex $P$ is involved. $\square$

\begin{Lemma} Suppose $(P,Q)$ is an edge of $\Gamma.$ Then 
$$d(\Gamma) = d_{PQ}-d_{P,Q}$$
\end{Lemma}

{\bf Proof.} This follows from Lemma 4.1 by breaking up the terms in the sum for $d(\Gamma)$ into two  classes: $(P,Q) \notin S$ and $(P,Q) \in S.$ $\square$

\begin{Lemma} Suppose $\Gamma'$ is obtained from $\Gamma$ by decreasing the weight of $P$ by $1.$ Then
$$d(\Gamma')=d(\Gamma) + d_P$$ 
\end{Lemma}
{\bf Proof.} Expanding and contracting by the vertex $P$, 
$$d(\Gamma')=(-a_P+1)\cdot d_P  - \sum \limits_{(P,Q) \in E(\Gamma)} d_{P,Q}=$$
$$=\Big( -a_P\cdot d_P  - \sum \limits_{(P,Q) \in E(\Gamma)} d_{P,Q}\Big) +d_P =d(\Gamma) + d_P $$
$\square$

\vskip 0.5cm

Generalizing the behavior of the weighted graphs of exceptional curves under blow-ups, we define the blow-up operations on the set of all weighted forests as follows.  The blow-down is defined as the operation opposite to blow-up. 

\begin{Definition}
1) Suppose $P$ is a vertex of $\Gamma.$ Then the blow-up of $\Gamma$ at $P$ is a forest with the set of vertices $V(\Gamma)\cup \{R\}$, the set of edges $E(\Gamma) \cup \{(PR)\}$;  $a_R=-1$, $a_P$ is decreased by $1$, all other weights do not change.

2) Suppose $(P,Q)$ is an edge of $\Gamma.$ Then the blow-up of $\Gamma$ at $P$ is a forest with the set of vertices $V(\Gamma)\cup \{R\}$, the set of edges $E(\Gamma) \cup \{(P,R), (R,Q)\}\setminus \{(P,Q)\}$; $a_R=-1,$ $a_P$ and $a_Q$ are decreased by $1$, all other weights do not change.
\end{Definition}

\begin{Lemma} The determinant of the graph is preserved by blow-ups and blow-downs.
\end{Lemma}

{\bf Proof.} Obviously, it is enough to prove this for the blow-ups.

Case 1. Suppose $P\in V(\Gamma)$ is blown up, and the new vertex is $R$. Suppose the new graph is $\Gamma'$. Expanding by vertex $R$, we get

$$d(\Gamma')=d(\Gamma'_R)-d(\Gamma'_{PR}) =d(\Gamma'_R)- d(\Gamma_P)$$
By the weight increment formula, 
$$d(\Gamma'_R)=d(\Gamma) + d(\Gamma_P),$$
which implies the result.

Case 2. Suppose $(P,Q)   \in E(\Gamma) $ is blown up, and the new vertex is $R.$  Suppose the new graph is $\Gamma'$.  Because $\Gamma$ is a forest, $\Gamma_{PQ}$ is a disjoint union of two forests. We will denote the one that contains $P$ by $P\Gamma$ and the one that contains $Q$ by $Q\Gamma.$ Similarly, $\Gamma'_R$ is a disjoint union of  $P\Gamma'$ and $Q\Gamma'.$

Expanding by the vertex $R,$ we get
$$d(\Gamma')=d(\Gamma'_R)-d(\Gamma'_{P,R})- d(\Gamma'_{Q,R})$$

By the multiplicativity property,
$$d(\Gamma'_R)=d(P\Gamma') \cdot d(Q\Gamma'),$$ 
$$d(\Gamma'_{P,R}) = d(P\Gamma'_P) \cdot d(Q\Gamma'),$$ 
$$d(\Gamma'_{Q,R}) = d(Q\Gamma'_Q) \cdot d(P\Gamma').$$  

Note that $P\Gamma'_P = P\Gamma_P.$ By the weight increment formula,
$$d(P\Gamma')=d(P\Gamma)+d(P\Gamma_P).$$

Likewise, $d(Q\Gamma')=d(Q\Gamma)+d(Q\Gamma_Q).$ Putting this together, we get
$$d(\Gamma')=(d(P\Gamma) + d(P\Gamma_P))(d(Q\Gamma)+ d(Q\Gamma_Q))-d(P\Gamma_P)(d(Q\Gamma)+ d(Q\Gamma_Q))- $$
$$-d(Q\Gamma_Q)(d(P\Gamma) + d(P\Gamma_P))=d(P\Gamma)d(Q\Gamma) - d(P\Gamma_P)d(Q\Gamma_Q)$$

On the other hand, expanding $d(\Gamma)$ by the edge $(P,Q),$ we get
$$d(\Gamma) = d(\Gamma_{PQ})- d(\Gamma_{P,Q}).$$
By the multiplicativity property, this equals the above  expression for $d(\Gamma')$, which completes the proof of the lemma. $\square$

\begin{Corollary} The determinants of vertices and edges are preserved by blow-ups and blow-downs that do not destroy them.
\end{Corollary}
{\bf Proof.} Most blow-ups and blow-downs, that do not destroy the corresponding vertex or edge, do not involve it. For those that do, we note that blowing up an edge $(P,Q)$ and then removing $P$ is the same as removing $P$ and then blowing up $Q$. Also, when a vertex $P$ is blown-up and then removed, the determinant of the resulting graph is the same as the determinant of the original graph with vertex $P$ removed, because the new vertex has weight $-1$. $\square$

\vskip 0.5cm
The above corollary means that just like the coefficients of the $\bar{K}$ labels, the determinant labels on vertices and edges of the graph of exceptional curves on $Z$ do not change once the point or edge is created (the edge, however, may get destroyed). While we are ultimately more interested in the labels of vertices, we need to keep track of the edges as well. The following lemma is very significant, as it describes the determinants of  vertices and edges formed after  the blow-up in terms of the nearby determinants and the determinant of the entire graph. The existence of these formulas is not intuitively obvious, and the formulas themselves are very important for applications.

\begin{Lemma} Suppose $\Gamma'$ is obtained from $\Gamma$ by a single blow-up. Suppose $d(\Gamma')=d(\Gamma)=d.$

Case 1. Suppose $P\in V(\Gamma)$ is blown up, and the new vertex is $R$. Then $d_R=d_P+d$ and $d_{PR}=d_P+d.$ Here $d_R$ and $d_{PR}$ are the determinants of a vertex and an edge in $\Gamma',$ and $d_P$ is the  determinant in $\Gamma$ or $\Gamma'$ (which are the same by the previous lemma).

Case 2. Suppose $(P,Q)   \in E(\Gamma) $ is blown up, and the new vertex is $R.$ Then 
$$d_R=2d_{PQ}+d_P+d_Q-d,$$
$$d_{PR}=d_P+d_{PQ},$$
$$d_{QR}=d_Q+d_{PQ}.$$
\end{Lemma}

{\bf Proof.} In Case 1, one immediately sees that $d_{PR}=d_R,$ because the weight of $R$ is $(-1).$ 
Then the weight increment formula implies that $d_R= d +d_P .$

In Case 2, using the notation and work from the previous lemma, we have
$$d_R=d(\Gamma'_R) =d(P\Gamma')d(Q\Gamma')=(d(P\Gamma)+d(P\Gamma_P))\cdot (d(Q\Gamma)+d(Q\Gamma_Q)) =$$
$$d_{PQ}+d_P+d_Q+d(P\Gamma_P)d(Q\Gamma_Q)=d_{PQ}+d_P+d_Q+(d_{PQ}-d).$$
This proves the formula for the vertex. For the edges, it is enough to prove the formula for $d_{PR}$.
$$d_{PR}=d(P\Gamma')\cdot d(Q\Gamma)= (d(P\Gamma)+ d(P\Gamma_P))d(Q\Gamma)= d_{PQ}+d_P. \ \square$$

\vskip 0.5cm

In what follows, the weighted forest $\Gamma$ will be the tree of exceptional curves on  $Z,$ which is a resolution of singularities of a counterexample to the Jacobian Conjecture with $X=P^2.$ We define the determinant label of an exceptional curve $P$ as $d_P.$ To help distinguish the different labeling, we will call the coefficient of $P$ in $\bar{K}_{Z}$ the $\bar{K}-$label of $P.$

The relevance of the determinant labeling stems from the following theorem.

\begin{Theorem} 
Suppose $T$ is any collection of type 1 curves, and  $\Gamma_T$ is the subgraph obtained from $\Gamma$ by removing the curves from $T$ and all affected edges. Then $d(\Gamma_T)<0.$ In particular, $d_P<0$ for any type 1 curve $P.$ 
\end{Theorem}

{\bf Proof.} Because of the Hodge Index Theorem, if $d(\Gamma_T)\geq 0,$ the restriction of the intersection form to the span of the classes of exceptional curves not in $T$ is negative semi-definite. On the other hand,   the Ample Ramification Theorem implies that some linear combination of curves of type $2,3,4$ on $Z$ has positive self-intersection.  $\square$

\begin{Theorem} Every $\bar{K}-$negative curve $P$ with negative $d_P$ has a $\bar{K}-0$ curve as its ancestor.
\end{Theorem} 

{\bf Proof.} For a curve $P,$ we will denote its $\bar{K}$ label by $b_P.$ We will identify the curves with vertices on the graphs. Suppose we perform any number of blow-ups without creating a $\bar{K}-0$ curve. This means that $b_P\leq -1$ for all curves $P$ involved. We will prove by induction on the number of blow-ups the following two conditions.

1) For all vertices $P$, $d_P+b_P\geq -1.$ Note that this implies that $d_P\geq 0.$

2) For all edges $PQ,$ $d_{PQ}\geq 0.$

The base of induction is easy to check. To prove the step, we need to consider two cases of blow-ups.
Note that $d(\Gamma)=d=-1.$

Case 1. Blowing up a vertex $P$ to get a new vertex $R$.
$$d_R+b_R=(d_P-1)+(b_P+1)=d_P+b_P\geq -1; \ \ d_{PR}=d_R\geq 0.$$

Case 2. Blowing up an edge $PQ$ to get  a new vertex $R.$
$$d_R+b_R=2d_{PQ}+d_P+d_Q+1+b_P+b_Q=$$
$$=(d_P+b_P)+(d_Q+b_Q)+1+2d_{PQ} \geq (-1) + (-1)+1+0=-1. $$
For the edges,
$$d_{PR}=d_P+d_{PQ}\geq 0+0=0.$$

 $\square$

\begin{Example} Starting from $X=P^2,$ we blow up a point at infinity, and then a point on the new curve. The resulting surface has the following graph of exceptional curves.

$$1\ \ \ \ _0\ \ \ \  0\ \ \ _{-1}-1$$
\vskip -0.9cm
$$\circ \!\! - \!\!\! - \!\!\! - \!\!\! - \!\!\! - \!\!\! - \!\! \circ \!\! - \!\!\! -\!\!\! - \!\!\! -\!\!\! - \!\!\! - \!\! \circ$$
\vskip -0.8cm
$$-2\ \ \ \ \  -1 \ \ \ \  \ \ \ \  0\ \ $$
\vskip -0.7cm
$$\  0\ \ \ \ \  \ -2 \ \ \ \ \ -1$$
Here the labels above the graph are the edge and vertex determinants, the labels right below the graph are the $\bar{K}-$labels, and the labels on the bottom line are the self-intersection labels. If one blows up the point of intersection of the last two blown-up curves, one gets a curve with negative $\bar{K}-$label and negative determinant label.  Other such curves could be obtained by further blowups.          
\end{Example}

The following theorem shows that, up to a polynomial automorphism of the original $A^2\subset X,$ all type $1$ curves are obtained as in the example above.

\begin{Theorem} Suppose $Z$ is a resolution of singularities of a counterexample to the Jacobian Conjecture. Suppose $P$ is $\bar{K}-$negative curve with negative $d_P$. Suppose $E$ is the first $\bar{K}-0$ curve that has to be created to create $P$ (the oldest $\bar{K}-0$ ancestor of $P$). Then there exists a sequence of blow-ups and blow-downs at infinity, that transform $Z$ into $Z',$ which has the same graph as the example above, with $E$ being the $\bar{K}-0$ curve. (Note that the curve with the $\bar{K}-$label $(-2)$ may differ from the original  pullback of infinity on $Z$).
\end{Theorem}

{\bf Proof.} We will follow the way $P$ was created, changing the pullback of infinity as needed, to decrease the number of ancestors of $E$. The key idea is that one can rule out the possibility of creating $P$ in some choices of blow-ups.

We start with blowing up a point at infinity. If then the new point is blown-up on the last curve, we are in the situation of the example above, and the result is obviously true. So the only non-trivial option left is to blow the point of intersection of the two curves, creating the following graph:
$$1\ \ \ \ _1\ \ \ \  2\ \ \ \ _{0}\ \ \ \ 0$$
\vskip -0.9cm
$$\circ \!\! - \!\!\! - \!\!\! - \!\!\! - \!\!\! - \!\!\! - \!\! \circ \!\! - \!\!\! -\!\!\! - \!\!\! -\!\!\! - \!\!\! - \!\! \circ$$
\vskip -0.8cm
$$-2\ \ \ \ \  \ -3 \ \ \ \ \ \  -1\ $$
\vskip -0.7cm
$$ -1\ \ \ \ \ \ -1 \ \ \ \ \  \ -2\  $$

The next step in the sequence of creation of $E$ must either blowing up the left edge, the middle point, or the right edge of the above graph.

If the left edge is blown up, we end up in the part of the graph with all positive edge labels, and  the inequality $d_Q+b_Q\geq 0$ satisfied for all the vertices (except the one that used to be the middle vertex). The induction argument, similar to the one in the proof of the last theorem,  implies that even if we create a $\bar{K}-0$ curve, its determinant label, and the label of the edge created with it, will be non-negative. Since we have to go back to the ``negative $\bar{K}$" territory, we will never get a curve with negative $\bar{K}$ and determinant labels. 

If the middle point is blown (which corresponds to blowing up a new point on the last blown-up curve) then the newly blown curve can serve as a new pull-back of infinity, which is closer to $E$.

So the only case left is when the right edge is blown up. After that we again have three possibilities. Blowing up the left edge can be discarded the same way as above. If we blow up the last vertex, we will be forced to blow up the newly blown-up vertex, creating a new curve that can serve as a pull-back of infinity.  

In general, suppose the creation sequence for $E$ starts with blowing up the right edge exactly $k$ times. This leads to the following graph:
$$1\ \  \ \ _1\ \ \ \  2\ \ \ _2\hskip 1.6cm  (k+2) \ \ _{0}\ \ 0$$
\vskip -0.88cm
$$\circ \!\! - \!\!\! - \!\!\! - \!\!\! - \!\!\! - \!\!\! - \!\! \circ \!\! - \!\!\! -\!\!\! - \ \ \ \  ... \ \ \ \ -\!\!\! -\!\!\!  \circ  \!\!\! -\!\!\! - \!\!\! - \!\!\! -\!\!\! -\!\!\!- \!\!\circ$$
\vskip -0.8cm
$$\ \ -2\ \ \  \ -3 \hskip 1.8cm  -(k+3) -1\ \  $$
\vskip -0.7cm
$$\ \  -1\ \ \ \ -2 \hskip 2cm
-1\ \ -(k+2)$$

After this we are forced to blow up the second from the right vertex, and then blow up the resulting vertex, and so on, until we get a curve with the $\bar{K}-$label $(-2)$. Notice that this curve $Q$ can serve as a new pull-back of infinity, because all other curves can be contracted one after another to a smooth surface, which has to be $P^2,$ since the anticanonical class is $-3Q$ and $Q^2=1.$ Note that this $Q$ is closer to $E$ that the original pull-back of infinity.  $\square$

We are now going to extend the above results to several $\bar{K}-$negative curves, with the condition that removing them all from the graph of exceptional curves produces a graph with negative determinant.

First, we need to introduce some new notation. Suppose $X=P^2$ is the standard compactification of $A^2$, $Z$ is obtained from $X$ by a sequence of blowups outside $A^2$; $\Gamma$ is the (dual) graph of the exceptional curves on $Z,$ $E=\pi^{-1}(\infty)$, where $\pi :Z \to X$ is the natural map. Suppose $P$ is an exceptional curve on $Z$, and a vertex of $\Gamma.$ We denote by $d'_P$ the determinant of $\Gamma_{E,P}$. We denote by $u_P$ the coefficient of $P$ in the formula that expresses $\pi^*(\infty)$ in terms of the exceptional curves on $Z$: $\pi^*(\infty) = E + \sum_{P\neq E} u_{P}P$. Note that $u_P$ is a natural number.

\begin{Lemma} In the above notation, $u_P^2=d_P+d'_P$.
\end{Lemma}

{\bf Proof.} While an inductive argument is possible, the following direct proof is easier. Define
$L=E+...+U_PP+..,$ the class of the pullback of $(\infty)$ on $Z.$ Note that $L$ intersects by $1$ with $E$ and by $0$ with all other exceptional curves. If we change a basis of the lattice of the divisor classes on $Z$, replacing $P$ by $L$, the Gram matrix of the minus-intersection form will be multiplied by $u_P^2$. Thus
$$(-1)\cdot u_P^2= \det \left( \begin{array}{ccc}
(-E^2)&... &-1\\
...&\ &\ \\
-1&\  &-1
\end{array} \right), $$ 
where the last row and column of the matrix on the right correspond to $L$ and the first correspond to $E$. Note that all the coefficients in the last row and column are zero unless indicated , and that removing these row and column produces the matrix of $\Gamma_P$. Using cofactor expansion w.r.t. the last row, and then (for one of the terms) the first row, we get
$-u_P^2=-d_P-d'_P$. $\square$

\begin{Lemma} Suppose $P$ and $Q$ are on the different sides of $E$ (i.e. $E$ belongs to the pass in the tree $\Gamma$ that connects $P$ and $Q$). Then the determinant $d_{P,Q}$ of the graph $\Gamma_{P,Q}$ is given by the following formula:
$$d_{P,Q} = u_Pu_Q-d_Pd_Q$$
\end{Lemma}

{\bf Proof.} We will denote $d_{P,Q}$ by $_Pd_Q.$ The notation $_{C_1C_2}d_Q$ will stand for the determinant of the graph obtained by removing the vertex $Q$ ad the edge $C_1C_2$.

We first perform a sequence of blowups to create $Q$, without producing any ancestors of $P.$ Then we perform a sequence of blowups that creates $P,$ keeping track of the determinants of vertices and edges of the new graph with the point $Q$ removed.  Note the following:

1) $_Ed_Q=d'_Q$;

2) $d(\Gamma_Q)=d_Q$;

3) For all ancestors $C$ of $P$ and all edges $C_1C_2$ involved in creating $P,$ their determinant labels $_Cd_Q$ and $_{C_1C_2}d_Q$ are linear combinations of $d'_Q$ and $d_Q.$ Let us denote the coefficients of these linear combinations by $a_C,$ $b_C,$ $a_{C_1C_2},$ $b_{C_1C_2}$ so that
$$_Cd_Q=a_C\cdot d'_Q+b_C\cdot d_Q,\ \ \ \ \ \ \ \  _{C_1C_2}d_Q=a_{C_1C_2}\cdot d'_Q+b_{C_1C_2}\cdot d_Q $$

4) From the formulas for the determinants of blowups (Lemma 5.7) the coefficients $a_C$ and $a_{C_1C_2}$ behave like $d_C$ and $d_{C_1C_2}$ for a graph with the determinant $d=0$, and initial value $1$. By induction, one can easily prove that $a_C=d_C+d'_C,$ which equals $u_C^2.$

5) From the formulas for the determinants of blowups (Lemma 5.7) the coefficients $b_C$ and $b_{C_1C_2}$ behave like $d_C$ and $d_{C_1C_2}$ for a graph with the determinant $d=1$, and initial value $0$. By induction, one can easily prove that $b_C=d'_C.$

6) So, $_Cd_Q=u_C^2\cdot d'_Q +d'_Q=u_C^2(u_Q^2-d_Q)+ (u_C^2-d_C)d_Q = u_C^2u_Q^2-d_Cd_Q$ $\square$

\begin{Lemma} In the above notation, for a curve $P$ with $\bar{K}$ label $\bar{k}_P <0$ and $d_P<0$ the following inequality is satisfied:
$$u_P^2\geq |d_P|$$
\end{Lemma}
{\bf Proof.} Because $d_P<0,$ we need to prove that $u_P^2+d_P\geq 0.$ Note that $u_P^2+d_P=2d_P+d'_P.$ For all ancestors $C$ of $P$, define $l_C=2d_C+d'_C$ and $l_{C_1C_2} = 2d_{C_1C_2} + d'_{C_1C_2}$.  The labels $l_C$ and $l_{C_1C_2}$ behaves like those in the graph with $d=-1$. Like in the proof of Theorem 5.2., one can prove by induction the following two statements:

1) $l_C+\bar{k}_C \geq -1$

2) $l_{C_1C_2} \geq 0$

So $l_P+\bar{k}_P \geq -1$. Because $\bar{k}_P <0,$ $l_P \geq 0.$ $\square$

\begin{Theorem}  Suppose $P$ and $Q$ have negative $\bar{k}$ labels and determinant labels and lie  on the different sides of $E$. Then $d_{P,Q} \geq 0.$
\end{Theorem}  
{\bf Proof.} This follows directly from Lemma 5.9 and Lemma 5.10. $\square$

\begin{Corollary} Suppose $\{E_1,...,E_k\}$ is any collection of exceptional curves of $Z$ with negative $\bar{k}$ labels, such that removing them all from the graph of exceptional curves produces a graph with negative determinant. Then they all lie on the one side of $E$.
\end{Corollary}

\begin{Remark}
The above Corollary can be applied to the set of all curves of type $1$, as defined in sections 1--3. As a result, we conclude that all of these curves are on one side of $\pi^{-1}(\infty)$.
\end{Remark}

\section{Changing the Target Surface}

In order to better understand the structure of a possible counterexample to the Jacobian Conjecture, one should not be content with just having  the $P^2$ as the compactification of the target plane. Indeed, we proved in  Section 3 that for the compactification $Y=P^2$ all curves that are contracted to a point at infinity (including the pull-back of the line at infinity) are mapped to one point. So all the ``action" is in (or, rather, above) this point, and we should blow it up to better understand what is going on. In fact, it makes sense to keep blowing up the image of the pull-back of the line at infinity until this image becomes a curve.

In what follows, we will use the following conventions: $X$ is compactification of the origin $A^2$ ($X=P^2$, unless specified otherwise), $Y$ is the compactification of the target $A^2,$ $Z$ is the resolution of the corresponding map, $Z\to W \to Y$ is the Stein factorization.   The maps are $\pi: Z\to X$, $\phi: Z \to Y,$ $\tau: Z\to W,$ and $\rho: W\to Y.$ When we consider several compactifications at the same time, these varieties and maps acquire matching indices (e.g. $\rho_3$ is the map from $W_3$ to $Y_3.$). The log-ramification divisor on $W$ is $\bar{R}=\sum r_iR_i,$ where $R_i$ are all di-critical curves and $r_i$ are the corresponding ramification indices.

\begin{Definition}
 For $n\geq 0$ we denote by ${\mathcal F}_n$ the smooth compactification of $A^2$ that contains exactly two exceptional curves, $F_0$ and $F_n,$  that are smooth and satisfy $F_0F_n=1$, $F_n^2=0$, $F_0^2=-n$.
\end{Definition}

These compactifications ${\mathcal F}_n$ appear naturally as the compactifications $Y$ of the target $A^2.$ For example, starting from $Y= P^2$ and  blowing up $\phi (\pi^{-1} (\infty)), $ one gets ${\mathcal F}_1$. The curve $F_1$ is the strict pullback of the line at infinity on $Y$ and the curve $F_0$ is the blown-up curve. One can show (see below) that on this new compactification, that we will call $Y_1,$ the map $\phi_1$ sends $\pi_1^{-1}(\infty)$ to the point of intersection of $F_0$ and $F_1$. (Here $\pi_1$ is the map from $Z_1$ to $X_1=X=P^2$). Easy calculations show that the $\bar{K}$-labels of $F_0$ and $F_n$ are $-1$ and $-(n+1)$ respectively.

\begin{Theorem} Suppose we have any compactification of a counterexample to the Jacobian Conjecture with $X=P^2$ and $Y={\mathcal F}_n$. Then $\phi(\pi^{-1}(\infty))$ must be a point on $F_n$ (it may be the point of intersection of $F_0$ and $F_n$). Moreover, all di-critical curves on $W$ pass through $\tau(\pi^{-1}(\infty))$.
\end{Theorem}

{\bf Proof.}  If $\rho (\pi^{-1}(\infty))$ is $F_0$ or a point on $F_0,$ then $\phi^* (F_n)$ is a linear combination of exceptional curves on $Z$ that does not include $\pi^{-1}(\infty).$ This is impossible because $(\phi^* (F_n))^2 = (\deg (\rho)) \cdot F_n^2=0.$

Recall that we have at least one di-critical curve. As in Theorem 3.2, consider the intersection of a di-critical curve $R_k$ on $W$with the log-canonical class. We use the adjunction formula for $\rho $ and the adjunction inequality for $R_k$ on $W$ from Section 4. Here $\bar{R} = \sum_i r_iR_i$ is the di-critical ramification divisor on $W,$ $\rho(R_k)$ is the reduced image of $R_k.$
$$-1\leq R_k \cdot \bar{K}_W =f_k \rho(R_k) \cdot  \bar{K}_Y + R_k\cdot \bar{R} $$

Note that $\rho(R_k) $ must intersect at least one of the curves $F_0$ and $F_n$, so $ \rho(R_k) \cdot  \bar{K}_Y \leq -1$.   So $R_k \cdot \bar{R} \geq 0$. Thus for every point $w$ on $W$ that is an intersection of a di-critical curve with the other exceptional curves, $(\sum \limits_{i | w\in R_i}  r_iR_i)^2\geq 0,$ which implies that $\pi^{-1}(\infty) $ is among the irreducible components of $\tau^*(\sum \limits_{i | w\in R_i}  r_iR_i )$, thus $\tau(\pi^{-1}(\infty))=w.$ $\square$

Suppose we have a counterexample to the Jacobian conjecture, compactified and resolved, so that $X$ and $Y$ are isomorphic to $P^2.$ Blowing up the point $\phi (\pi^{-1} (\infty)), $ we get a compactification $Y_1={\mathcal F}_1.$ If the point $\phi_1 (\pi_1^{-1} (\infty)) $  is the intersection of the exceptional curves $F_0$ and $F_1$ on $Y_1$, we blow it up and contract the strict preimage of $F_1,$ obtaining the compactification $Y_2={\mathcal F}_2.$ If the point $\phi_2 (\pi_2^{-1} (\infty)) $  is the intersection of the exceptional curves $F_0$ and $F_2$ on $Y_2$, we again blow it up and blow down the strict pullback of the exceptional fiber to get $Y_3 ={\mathcal F}_3,$ and so on. We will continue to do this until for some $n$ we get $\phi_n(\pi_n^{-1}(\infty))$ to be on $F_n$ and not on $F_0.$  

Then we blow up  $\phi_n(\pi_n^{-1}(\infty))$ on $Y_n.$ Note that $n\geq 2,$ so the $\bar{K}$-label of the newly blown-up curve is less than or equal to $-2.$ Thus the same argument as in the theorem above implies that the pullback of the line at infinity is mapped by $\tau$ to a point, and all di-critical curves pass through it. The image by $\rho$ of this point must be on the newly blown-up curve. If it is not on the strict pullback of $F_n,$ we can blow down the strict pullback of $F_n$, so that the new compactification is ${\mathcal F}_{n-1}$. We will keep doing this until $\phi(\pi^{-1}(\infty))$ is on the point of intersection of pullback of the exceptional fiber and the newly blown-up curve. This implies that before the last blowup the images of all di-critical components on $Y$ ``touched" the exceptional fiber. (The quotation marks are due to the fact that these images are uni-branched, but not necessarily smooth at $\phi(\pi^{-1}(\infty))$).

The following theorem determines the  structure of the variety $W_n$ and the map $\rho_n$ in the case when $\phi_n(\pi_n^{-1}(\infty))$ is a point on $F_n$ and not on $F_0$.

\begin{Theorem} Suppose $X=P^2,$ $Y={\mathcal F}_n,$ and   $\phi_n(\pi_n^{-1}(\infty))$ is a point on $F_n\setminus F_0.$ We denote by $R_k$ the di-critical curves on $W,$ by $E_j$ the preimages of $F_n$ and by $D_i$ the preimages of $F_0.$ Then the following are true.

1) All $E_j$ contain $w=\rho(\pi^{-1}(\infty))$.

2) Every $D_i$ intersects exactly one $E_j$ at a cyclic quotient singularity and is smooth elsewhere (generalized simple normal crossing). The map $\rho$ is degree one  on each $D_i$ (but it may be ramified along $D_i$).

3) All $R_k$ contain $w$.

4) There are no other points of intersections of exceptional curves on $W.$ 
\end{Theorem}

{\bf Proof.} Statement (3) was proven in Theorem 6.1. Statement (2) follows from  Theorem 4.6 and Theorem 4.7 (a).  To prove statement (1) note that if some $E_j$ does not contain $w,$ it cannot intersect with any $R_k$. So Theorems 4.6 and 4.7(a) apply to it. Because $E_j\rho^*(F_0) >0,$ it must intersect with some $D_i$. Since that $D_i$ can intersect no other exceptional curve, and nether can $E_j$, we get a contradiction with the connectedness of the set of exceptional curves on $W.$ Finally, statement (4) follows from the fact the set of exceptional curves on $W$ is an image of a tree of curves on $Z.$ $\square$

\begin{Remark} 1) Note that in the above theorem we have not used the fact that the images of the di-critical curves ``touch" $ F_n$.  Perhaps this can be explored by studying $W$ for the surface $Y$ obtained by blowing up $\rho(w)$.

2) One can calculate $E_j\cdot F_i$ for intersecting $E_i$ and $F_j$ based on the ramification indices.

3) If $f_j=1,$ then $E_j$ must be smooth outside $w$ (obtained by applying Theorem 4.6 after blowing up $\phi(\pi^{-1}(\infty))$ several times); in this case there must be only one $F_i$ intersecting this $E_j$. For $f_j>1$ there may be several $F_i$ intersecting $E_j$ and $E_j$ may have singularities at the points that are mapped by $\rho$ to $\rho (w)$. 
\end{Remark}

We will now continue to blow up $\phi(\pi^{-1}(\infty))$ until it becomes a curve $F$ on $Y$. Note that because the $\bar{K}$-label of $E=\pi^{-1}(\infty)$ is $-2,$ the $\bar{K}$-label of $F$ can only be $-1$ or $-2,$ depending on whether the ramification index at $E$ is $2$ or $1$. Thus, no curves with positive $\bar{K}$-label are ancestors of $F.$ So $Y$ obtained by this procedure only contains curves with non-positive $\bar{K}$-labels. The following theorem indicates that it must contain $\bar{K}-0$ curves. Please note a notation change: the curves $E_i$ and $F_i$ now mean non-di-critical exceptional curves on $W$ and $Y$ respectively.

\begin{Theorem} Suppose $Y$, $E$ and $F$  are  the varieties described above, and $X=P^2,$ $Z,$ $W$ are the corresponding varieties, with the usual notation for the morphisms between them. Then for every di-critical curve $R_k$ on $W$ the following are true.

1) The point of intersection of $\rho(R_k)$ with $Y\setminus A^2$ belongs to exactly one exceptional curve $F_k.$

2) This curve $F_k$ has $\bar{K}$-label $0$.

3) The curve $F_k$ intersects exactly one other exceptional curve, with $\bar{K}$-label $-1.$
\end{Theorem}

{\bf Proof.} Suppose $\rho(R_k)$ intersects $Y\setminus A^2$ at a point $y_k,$ and $R_k$ intersect the union of other exceptional curves on $W$ at $w_k.$ Then if $y_k$ belong to at least one curve with negative $\bar{K}$-label, 
$$-1\leq R_k\cdot \bar{K}_W = R_k \cdot \rho^*(\bar{K}_Y) + R_k\cdot \bar{R} \leq -1 + R_k \cdot \bar{R} $$
Therefore $R_k\cdot \bar{R} \geq 0.$ The same is true for all di-critical curve passing through $w_k$,  so
$$(\sum \limits_{j|w_k\in R_j}   r_jR_j)^2 \geq 0$$

This implies that $(\sum \limits_{j|w_k\in R_j}  \tau^*( r_jR_j))^2\geq 0,$ which contradicts the fact that the determinant label of $E$ is positive. The same argument applies at any intermediate step to the curves $R_k$ that do not contain $\tau(\pi^{-1}(\infty))$.

When a $\bar{K}-0$ curve is created on $Y,$ it must be a result of blowing up a point on a curve with $\bar{K}$-label $-1.$ At this step, $\phi(\pi^{-1}(\infty))$ has to be the point of intersection of these two curves, otherwise we will get into the positive $\bar{K}$ territory on $Y$. So the statement of the theorem is true for $R_k$ at the intermediate step when $\rho(R_k)$ no longer contains $\phi(\pi^{-1}(\infty))$ ; it is easy to see that the subsequent blowups of $\phi(\pi^{-1}(\infty))$ do not change that. $\square$

\begin{Definition} For a hypothetical counterexample to the Jacobian Conjecture, we call the varieties $X,$ $Y,$ $Z,$ $W$ and morphisms between them, obtained by the sequence of blowups on $Y$ of $\phi(\pi^{-1}(\infty))$ until it becomes a curve, its resolution of the image of the line at infinity.
\end{Definition}

\begin{Remark}
The varieties in the resolution of the image of the line at infinity differ slightly from those described in Theorem 6.3, because we do not blow down any curves on $Y$, as we did during the modifications of ${\mathcal F}_n$. However, the proof and the conclusions of Theorem 6.3 remain true. Note also that the exceptional curves on $Y$ can be enumerated w.r.t. to the order of their creation. In particular, we can talk about the last $\bar{K}-0$ curve there. Note that the resulting $Y$ and $W$ are completely determined by the original pair of polynomials, while $Z$ may vary.   
\end{Remark}

 The above discussion allows to give an alternate proof of the Remark 5.1. Combined with the ``One Point at Infinity" theorem of Abhyankar (cf. \cite{AbhyankarTataLectures}) and the results of section 4, it gives the following result.

\begin{Theorem} Suppose a counterexample to the Jacobian Conjecture exists. Then there exists another counterexample, such that its resolution of the image of the line at infinity has the following properties, in addition to those described in Theorem 6.3.

1) If $F_l$ is the last $\bar{K}-0$ curve on $Y$, then  at least one image of a di-critical curve intersects it.

2) The curve $F=\phi(\pi^{-1}(\infty))$ on $Y$ is a curve of valency two. 

3) Removing $F$ from the graph of exceptional curves on $Y$ produces two subgraphs. One of them contains all the exceptional curves of type $1$ in the sense of sections 1--3. The other contains the curve $F_l$ and one of the following is true about it:

Case A) It consists of a single curve $F_l$; then $\phi(\pi^{-1}(\infty))$ is the $\bar{K}-(-1)$ curve it is connected to.

Case B) It consists of two curves: $F_l$ and the $\bar{K}-(-1)$ curve it is connected to; then  $\phi(\pi^{-1}(\infty))$ is  a $\bar{K}-(-2)$ curve, connected to the $\bar{K}-(-1)$ neighbor of $F_l$.
\end{Theorem}

{\bf Proof.} The ``One Point at Infinity" theorem of Abhyankar (cf. \cite{AbhyankarTataLectures}) implies that if for a given counterexample to the two-dimensional Jacobian Conjecture  the strict preimage of a generic line $L$ on the target $P^2$ is a curve on the origin $P^2$ with only one point at infinity, one can produce (precomposing with a suitable polynomial automorphism) a ``smaller" counterexample. In our terminology, the points of infinity are the points $\pi (\rho^{-1}(A_i)),$ where $A_i$ are points of intersection of $L$ and the images of the di-critical curves, and the intersection of $L$ with the line at infinity on the target $P^2.$ Thus for a given counterexample, we can precompose it with the polynomial automorphism to get another counterexample, where not all $\pi (\rho^{-1}(R_i))$ and $\pi (\rho ^{-1} (\infty))$ are the same point. This is equivalent to the condition that not all di-critical curves and curves of type $1$ in the terminology of sections 1--3  are on the same side of the curve $\pi^{-1}(\infty)$ on $X$. Note that this condition is preserved by the subsequent modifications of $X$. Note that by the Remark 5.1 all curves of type 1 must lie on one side of the curve $\pi^{-1}(\infty)$ on $X$. So at least one di-critical curve does not lie in the same connected component of the tree obtained by removing $\pi^{-1}(\infty)$.

Suppose $F'$ is the last created $\bar{K}-0$ curve on $Y$ that  intersects with the images of the di-critical curves. We claim that all points blown up after its creation were points of intersection of two exceptional curves on $Y$. Indeed, otherwise all images of the di-critical curves on $Y$ and the strict pull-back of the line at infinity of the original target plane would have been in the same connected component of the tree of the exceptional curves with $\phi(\pi^{-1}(\infty))$ removed. Because we stop the process as soon as $\phi(\pi^{-1}(\infty))$ is a curve, its valency is $1$ or $2$. If its valency is $1$, then Theorem 4.7 implies that $E=\tau(\pi^{-1}(\infty))$ is an end-curve in the graph of exceptional curves on $W$, which contradicts our assumption that not all di-critical curves and curves mapped to the original line at infinity are in the same connected component of the graph of exceptional curves on $X$. If the valency is $2$, then Theorem 4.7 implies that the valency of $E$ is two. Because all exceptional curves on $Y$ are ancestors of $F,$ the connected component of the graph of exceptional curves on $Y$ that does not contain the di-critical curves and the target line at infinity is a chain. Therefore Theorem 4.7 implies that every exceptional curve on $W$ that is mapped by $\rho$ to any of these curves has the same valency (one or two) as its image.  Therefore the map $\rho$ sends the connected component of the tree of curves on $W$ without $E$ which does not contain the curves of type $1$ to this chain of curves.  This again implies that all di-critical components and curves of type $1$ are on the same side of $E,$ a contradiction.

Therefore, after $F'$ was created we could only blow up points of intersections of two curves. This immediately implies that $F'=F_l.$ Recall that the $\bar{K}$ label of $F$ is $-1$ if the ramification at $E$ is $2$ and $-2$ if the ramification at $E$ is $1$. One can easily check by looking at the $\bar{K}$ labels on $Y$ that it is possible that a point of intersection of $F$ and its $\bar{K}-(-1)$ neighbor is blown up several times in the process, but no other points can be blown up until, possibly, one final blow-up that creates a $\bar{K}-(-2)$ curve. This implies the result. $\square$

\section{Informal Discussion and Observations}

First of all, by further blowups on $Y$ we can get a resolution on which the images of the  di-critical curves intersect the other exceptional curves transversally. Such resolutions were considered by Domrina and Orevkov (cf. \cite{Domrina}, \cite{DomrinaOrevkov}). Note that the maps of the splice diagrams that they considered have a nice geometric counterpart: if one contracts the chains of rational curves on $Y$ to contract all curves of valency two, the curves on $Y$ will have generalized simple normal crossings. Using the methods of section 4, one can easily show that the same will then be true for $W$, except for up to one cyclic quotient singularity on each of the di-critical curves. In a sense, our section 4 can be viewed as a strong argument for the importance of the splice diagrams.  It puts a geometric structure on top of the combinatorial structure of Domrina and Orevkov. In particular, one can try to find a counterexample to the Jacobian Conjecture by finding a suitable map between splice diagrams and then constructing the corresponding finite map of singular compactifications of $A^2.$  The restrictions obtained in section 6 are strong, but not obviously prohibitive, and they may tell us where to look for possible counterexamples.

It may be interesting to investigate the connections between our approach and other results on the Jacobian Conjecture (see, e.g. \cite{LeDungTrang}, \cite{Miyanishi}, \cite{Rabier2}). It is also clear that we have barely tapped into the powerful methods of modern birational geometry of surfaces. In particular, some ideas of the log Sarkisov Program (cf. \cite{BrunoMatsuki}) may help to modify the possible counterexamples. It should also be noted that the theory of rational curves on rational surfaces, in particular the methods and results of the paper \cite{KeelMcKernan}  may be helpful. Specifically, one can show that on $Y$ there cannot exist a covering family of rational curves that intersect the union of the curves at infinity and the images of the di-critical curves at only one point. If $D=\sum F_i$ is the support of $\rho (\bar{R})$, and $\rho^*D=\bar{R} +G$ (i.e. $G$ is the union of the co-exceptional curves on $W$) then on $W$ we have $\bar{K}_W + G= \rho^*(\bar{K}_Y+D)$. Therefore some multiple of  $\bar{K}_W + G$ must have a global section, an interesting condition to explore.

Finally, the following result seems to be true, but  will appear elsewhere, once the details are clarified and confirmed.

\begin{Theorem} Suppose $a$ and $b$ are any two fixed integers. Then there are only finitely many, up to the plane automorphisms, minimal graphs of exceptional curves that produce divisorial valuations with the center at infinity such that its determinant label is $a$ and its $\bar{K}$ label is $b$. In other words, such valuations are bounded, up to the plane automorphisms.
\end{Theorem}

\end{document}